\documentclass{jktr}
\usepackage{amstex}
\newcommand{\bdy}{\partial M}
\newtheorem{definition}{Definition}
\newtheorem{theorem}{Theorem}
\newtheorem{lemma}{Lemma} 
\input{epsf}

\title{A Mayer-Vietoris Theorem for \\ the Kauffman Bracket Skein Module}
\author{Walter F. LoFaro}
\address{Department of Mathematics,\\
         The University of Iowa,\\
         Iowa City, IA 52242, USA\\
         email:  lofaro@@math.uiowa.edu}
         
\makeatother         

\begin{document}

\maketitle

\begin{abstract}
The nth relative Kauffman bracket skein modules are defined and two theorems
are given relating them to the Kauffman bracket skein module of a 3-manifold.
The first theorem covers the case when the 3-manifold is split along a
separating closed orientable surface and the second theorem addresses the
case when the surface is nonseparating.
\end{abstract}

\section{Introduction}
The study of knots and links was invigorated in 1984 by the introduction
of the Jones polynomial \cite{vj}.  It allowed solution of many long standing
problems in knot theory including the Tate conjectures.  One drawback was
that the Jones polynomial is hard to compute for complicated knots and
links.  Then in 1987 Kauffman introduced the Kauffman bracket polynomial \cite{k}.
Easier to compute than the Jones polynomial, it carries essentially the
same information (modulo framing and orientation).  The fact that the
Kauffman bracket takes framings into account makes it an ideal tool for 
studying links in arbitrary 3-manifolds, which need not have a canonical
framing.

A major question at this time was how to consturct something similar to the
Jones polynomial for links in any 3-manifold.  One approach was via
topological quantum field theory.  Another appraoch was through skein
modules, introduced by Przytycki in 1989. \cite{prz}  The importance
of skein module theory was demonstrated when it was unified with topological
quantum field theory by Blanchet-Habiger-Masbaum-Vogel \cite{bhmv}, Lickorish 
\cite{l}, and Kauffman-Lins \cite{kl}.

In this paper we study the Kauffman bracket skein module $K(X)$ of an
oriented 3-manifold $X$.  $K(X)$ has been computed for few $X$.  In a
paper by Hoste and Przytycki a method for computing $K(X)$ based on a
handle decomposition of $X$ was introduced.  They found $K(X)$ when $X$
is a lens space \cite{hp1} and when $X=S^2 \times S^1$ \cite{hp2}.  Bullock
used this approach to determine $K(X)$
when $X$ is the complement of a $(2,2p+1)$-torus knot \cite{db}.  Unfortunately
this method seems too complicated when $X$ has high Heegaard genus.

Until now , one thing that has been missing fron the theory of Kauffman
bracket skein modules is a decomposition theorem in the spirit of the
Mayer-Vietoris theorem from homology theory.  Herein this problem is addressed.
Let $M$ be a smooth oriented 
3-manifold and $F$ a subsurface of $\bdy$.  Informally the nth
relative Kauffman bracket skein module of $M$, $K_n(M)$, is the Kauffman
bracket skein module of $M$ where the links have $n$ arc components ending
in $2n$ prescribed arcs in $F$.  Suppose we have two manifolds like $M$,
say $X_1$ and $X_2$, and $X=X_1 \cup_F X_2$.  We define a homomorphism

             \[\Phi:\bigoplus_{n=0}^{\infty}(K_n(X_1) \otimes K_n(X_2))
               \rightarrow K(X)\]

Intuitively, the map $\Phi$ is obtained by glueing an element of $K_n(X_1)$
to one of $K_n(X_2)$ to get an element of $K(X)$.  It is shown that $\Phi$
is onto and its kernel is described simply and topologically.

Suppose $L$ is a framed link in $X$ intersecting $F$ transversely in $2n$
arcs.  We can view $L$ as $L_1 \otimes L_2$ where $L_i = L \cap X_i$, $i=1,2$.
Suppose that near $F$ $L_i$ is vertical in $F \times  I$.  If $\sigma$ is a
$2n$ component framed braid over $F$ we can define $\sigma \cdot L_i$ as the
result of replacing $L_i \cap (F \times I)$ with $\sigma$.  The isotopy that
takes $L_1 \otimes L_2$ to $\sigma \cdot L_1 \otimes \sigma^{-1} \cdot L_2$
is called a braiding move.  Suppose some component of $L_i$ bounds a bigon
with one boundary component in $L_i$ and the other in $F$.  The isotopy
pulling that part of $L_i$ to the other side is a bigon move.  The kernel is 
the submodule generated by all expressions of the form $L-L^{\prime}$ where
$L^{\prime}$ is obtained from $L$ by a sequence of braiding and bigon moves.

\section{Definitions}
Suppose that $X$ is an orientable three-manifold and $F$
is a two-sided
orientable surface embedded in $X$. Let
$\psi:F \rightarrow F$
be a diffeomorphism and $X=X_1 \cup_{\psi} X_2$, where $X_i$ is also an orientable three-manifold 
for $i=1,2$. Of course, we could choose
coordinates so that $\psi$ is the identity; we assume this to be
the case throughout.  We will define the relative
Kauffman bracket skein modules for $X_i$, $K_n(X_i)$, and explore
the relationship between the $K_n(X_i)$ and the Kauffman bracket
skein module of $X$, $K(X)$.

Let $M$ be an orientable three manifold with boundary
$\partial M
=F$.  Let $a_{1}, \ldots ,a_{2n}$ be a collection of
disjoint embedded arcs on $F$.

\begin{definition}                                                
A framed arc in $M$ is a homogeneous embedding of the square
$I \times I$ into $M$ with $(I \times I) \cap \partial
M=(I \times \{0\}) \cup (I \times \{1\})$.
\end{definition}

We call the parts of the framed arc lying on $\partial M$ the
ends of the arc, and denote it by $\partial_0 L$. 

\begin{figure}[h]
\centering
\leavevmode
\epsfxsize = 2in
\epsfysize = 2in
\epsfbox{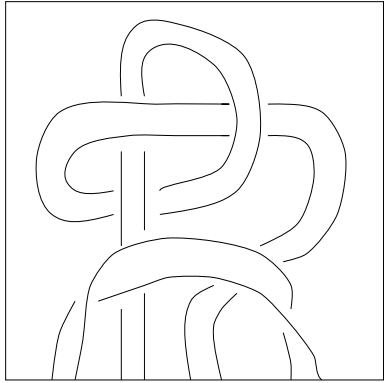}
\caption{Two framed arcs}
\end{figure}

Usually when we draw framed arcs and circles we just draw them as arcs.
In order to reflect the framing of the components, we use nugatory
crossings and we take the framing to be the blackboard framing.  See Figure 1 and Figure 2.
\begin{figure}[h]
\centering
\leavevmode
\epsfxsize = 2in
\epsfysize = 2in
\epsfbox{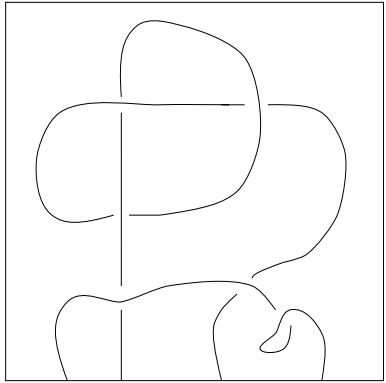}
\caption{Figure 1 represented as two arcs}
\end{figure}

\begin{definition}
A collection of $n$ framed arcs in $M$ is nicely embedded if
the ends of the framed arcs coincide with the $a_{j}$ in a 
one-to-one manner.
\end{definition}

We are now ready to define the relative Kauffman bracket skein
modules of the $X_{i}$, but let's first recall the definition of
the  Kauffman bracket skein module of $X$.  Let $R$ be the ring
of Laurent polynomials over the integers and let ${\cal L}(X)$ be
the set of all isotopy classes of links in $X$. The Kauffman
bracket skein module of $X$ is $R{\cal L}(X)/K$, where $K$ is the
submodule
generated by the standard skein and frame relations. 
See Figure 3.

\begin{figure}[h]
\centering
\leavevmode
\epsfxsize = 2in
\epsfysize = 2in
\epsfbox{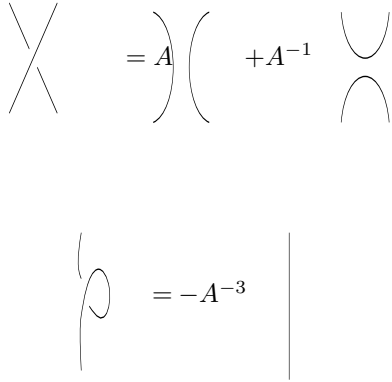}
\put(-90,30){=}
\put(-80,30){$-A^{-3}$}
\put(-100,120){$=A$}
\put(-55,120){$+A^{-1}$}
\caption{The skein and framing relations}
\end{figure}

For the
relative skein module we make the analogous definitions.  That
is, ${\cal L}_{n}(M)$ is the set of isotopy classes of nicely
embedded framed links with $n$ framed arcs and $K_{n}$ is the
submodule generated by the frame and skein relations with the
relations also applied to the framed arcs.  There is no
generality lost in restricting to nicely embedded framed arcs in
${\cal L}_{n}(M)$, for there is an isotopy of any framed arc
taking it to a nicely embedded one.  A proof of this would be
very much like the proof of Lemma 1.

\begin{definition}
The $n^{\textstyle{th}}$ relative Kauffman bracket skein module
of $M$, written $K_{n}(M)$ is $R{\cal L}(M)/K_{n}$. 
\end{definition}

Let $B_{2n}(F)$ be the $2n$-strand framed braid group over $F$. 
We want to define a group action of $B_{2n}(F)$ on ${\cal
L}_{n}(M)$.  Let
$L\in {\cal L}_{n}(M)$.  By taking a small product neighborhood
$\partial M \times I$ of $\partial M$ in $M$ it is possible
to isotope $L$ rel $\partial M$ so that $L \cap (\partial M
\times I)=(a_{1} \times I) \cup \cdots \cup (a_{2n} \times I)$. 
Thus there is no harm in assuming this about $L$ from the start.  

\begin{definition}
Let $\sigma \in B_{2n}(F)$. Define $\sigma \cdot L = (L-(L \cap
(\partial M \times I))) \cup \sigma$
\end{definition}

\begin{figure}[h]
\centering
\leavevmode
\epsfxsize = 3in
\epsfysize = 1in
\epsfbox{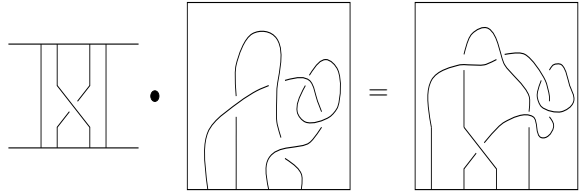}
\put(-80,35){$=$}
\caption{An example of the action of the $B_{2n}(F)$ on ${\cal L}_n(X_i)$}
\end{figure}

It is not hard to see that this really is a group action.  It is
clear that $(\sigma \tau) \cdot L = \sigma \cdot (\tau \cdot L)$,
and $e \cdot L = L$ by the assumption on $L$ immediately
preceding the definition.  Here $e$ is the identity braid.

Let $L \in {\cal L}_n(M)$ and suppose that
$H:(L \times I, \partial_0 L \times I) \rightarrow (M, \bdy)$ is an isotopy with $H_0(
\partial_0 L)=H_1(\partial_0 L) = \{a_0, \ldots , a_n\}$. This type of motion will be 
called a braiding move.
Lemma 2 will allow us to see that $H_1(L) = \sigma \cdot L$ up to isotopy rel $\bdy$ for
some $\sigma \in B_{2n}(F)$, where $L$ and $H_1(L)$ are considered as
elements of $K_n(M)$.

Now suppose that $L$ is a framed link in $X$ and $H$ is an
isotopy of $L$ in $X$ with $H_0(L) \cap F=H_1(L) \cap F$, and $H_t(\partial_0 L)=
H_t(L) \cap F$. Under these conditions $H_j(L) \cap X_i$ represents an
element of $K_n(X_i)$ for $j=0,1$.  Furthermore $H_1(L) \cap X_1
= \sigma \cdot (H_0(L) \cap X_1)$ and $H_1(L) \cap X_2 = \tau \cdot (H_0(L) \cap x_2)$ 
for some $\sigma ,\tau \in B_{2n}(F)$. By choosing the right coordinates on 
$\partial X_2$ we would actually have $\tau = \sigma^{-1}$.
Now let $U_n$ be the submodule of $K_n(X_1) \otimes K_n(X_2)$
generated by all elements of the form \[(a \otimes b) - (\sigma
\cdot a \otimes \sigma^{-1} b)\]
Finally set $U=\oplus_{n=0}^{\infty} U_n$.

\begin{definition}
Suppose that $L \in {\cal L}_n(X_i)$
has an arc component $c$ and there is a bigon bounded by $c$ and
some arc embedded in $F$. If $L$ does not intersect the interior
of this bigon
then an isotopy that pulls $c$ across to the other side of $X$
and fixes the other components of $L$ is called a bigon move.
\end{definition}

\begin{figure}[h]
\centering
\leavevmode
\epsfxsize = 3.5in
\epsfysize = 1in
\epsfbox{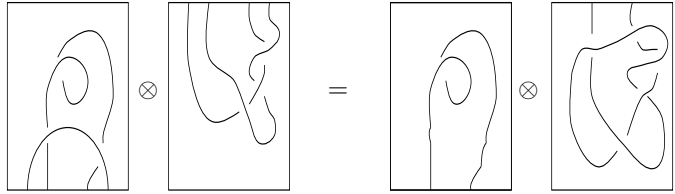}
\put(-203,36){$\otimes$}
\put(-59,36){$\otimes$}
\put(-131,36){$=$}
\caption{An example of a bigon move}
\end{figure}

Let $V$ be the submodule of $\oplus_{n=0}^{\infty}[K_n(X_1)
\otimes K_n(X_2)]$ generated by all elements of the form 
\[(a \otimes b) - (c \otimes d)\]
where $c \otimes d$ is the result of applying a bigon move to $a
\otimes b$.

\begin{theorem}
Let $X$ be an orientable three-manifold, $F \subset M$ a
closed orientable surface, and $\psi :F \rightarrow F$ a
diffeomorphism.  If $X=X_1 \cup_{\psi}
X_2$ then \[K(X) \cong  \frac{\oplus_{n=0}^{\infty}
K_n(X_1) \otimes K_n(X_2)}{U \vee V}  \]
\end{theorem}

\section{Lemmas}  
The first lemma we prove assures us that given any $L \in {\cal
L}(X)$ it is represented by a link that intersects $F$ in the
prescribed way.  Then we show that it is possible to
reparameterize $L$ so that, in the absence of bigon moves, an
isotopy of $L$ in $X$ "splits in a nice way" into isotopies in
$X_1$ and $X_2$.  Lastly we prove the theorem.

\begin{lemma}  Let $L \in {\cal L}_n(X)$.  There is an isotopy of
$L$ so that $L \cap X_i \in {\cal L}_n(X_i), i=1,2$ for some $n$.
\end{lemma}

\begin{proof}
Generically $L$ intersects $F$ in $2n$ arcs for
some $n$. Let $s_1, \ldots , s_{2n}$ be these arcs. Let $t_i$
be one point of $\partial a_i$ in $F$.  By the Isotopy Lemma
there is a diffeomorphism of $F$ taking
$s_i$ to $t_i$ for all $i$ and this diffeomorphism is isotopic to
the identity map on $F$. Call this isotopy $H:F \times I
\rightarrow F$. Next let
$b$ be a smooth bump function on $[-1,1]$ such that          
$b(-1)=b(1)=0$ and $b(0)=1$. We define another isotopy, called
$\hat H$ on a small product neighborhood $F \times [-1,1]$ of $F$
in $X$. Let $\hat{H}:(F \times [-1,1]) \times I \rightarrow    
[-1,1]$ be given by $\hat{H} (x,s,t)=H(x,b(s)t)$.  Finally it is
clear that after applying $\hat H$
to $L$ we can twist its band so that the altered $L$ intersects
$F$ in the $a_j$. 
\end{proof}

Now let's focus our attention on the characterization of the
submodules $U_n$.

\begin{lemma}
Let $M$ be an orientable three-manifold with 
boundary $F$.  Let $L \subset M$ be a properly embedded
(finite) collection of disjoint framed arcs and framed circles. 
Suppose there are $n>0$ framed arcs in $L$ and $\partial L =
\{a_1, \ldots , a_{2n}\}$, where the $a_j$ are as before. 
Suppose further that $H:L \times [0,1] \rightarrow M$ is an
isotopy of $L$ in
$M$ for which $H(L,1) \cap \partial M =\{a_1, \ldots , a_{2n}\}$. Then for some
$\sigma \in B_{2n}(F)$, $H(L,1)=\sigma \cdot L$ up to isotopy
rel
$\partial M$.
\end{lemma}
\begin{proof}
The strategy is to define a new isotopy $G:L \times [\delta,
\eta] \rightarrow M$ equivalent to $H$ in the sense that
$G(L,\eta)=H(L,\eta)$. $G$ will
also have the property that $G(L,\frac{\delta + \eta}{2})=\sigma
\cdot L$ for
$\sigma \in B_{2n}(F)$ and $G:L \times [\frac{\delta +
\eta}{2},1]
\rightarrow M$ is an isotopy rel $\partial M$.  As in Lemma 1, we
prove it for unframed nicely embedded links, but the result
carries through.

We want to choose a product neighborhood $\bdy \times [-1,1]$
with $\bdy = \bdy \times \{1\}$, so small that the arcs $L \cap
(\bdy \times [-1,1])$ have no critical points with respect to
height.  See Figure 6.  Stated more precisely, let $f_s:[- \epsilon, \epsilon]
\rightarrow \bdy \times [-1,1]$ be a parametrization of a
component of $(\bdy \times [-1,1])$ at time $s$ with $f_s(-
\epsilon) \in \bdy$.  Let $h:\bdy \times [-1,1] \rightarrow
[0,2]$ be given by $h(x,t)=1-t$.  Then we choose a neighborhood
so small that $h \circ f_s$ is strictly increasing for all $s$.

\begin{figure}
\centering
\leavevmode
\epsfxsize = 3in
\epsfysize = 2in
\epsfbox{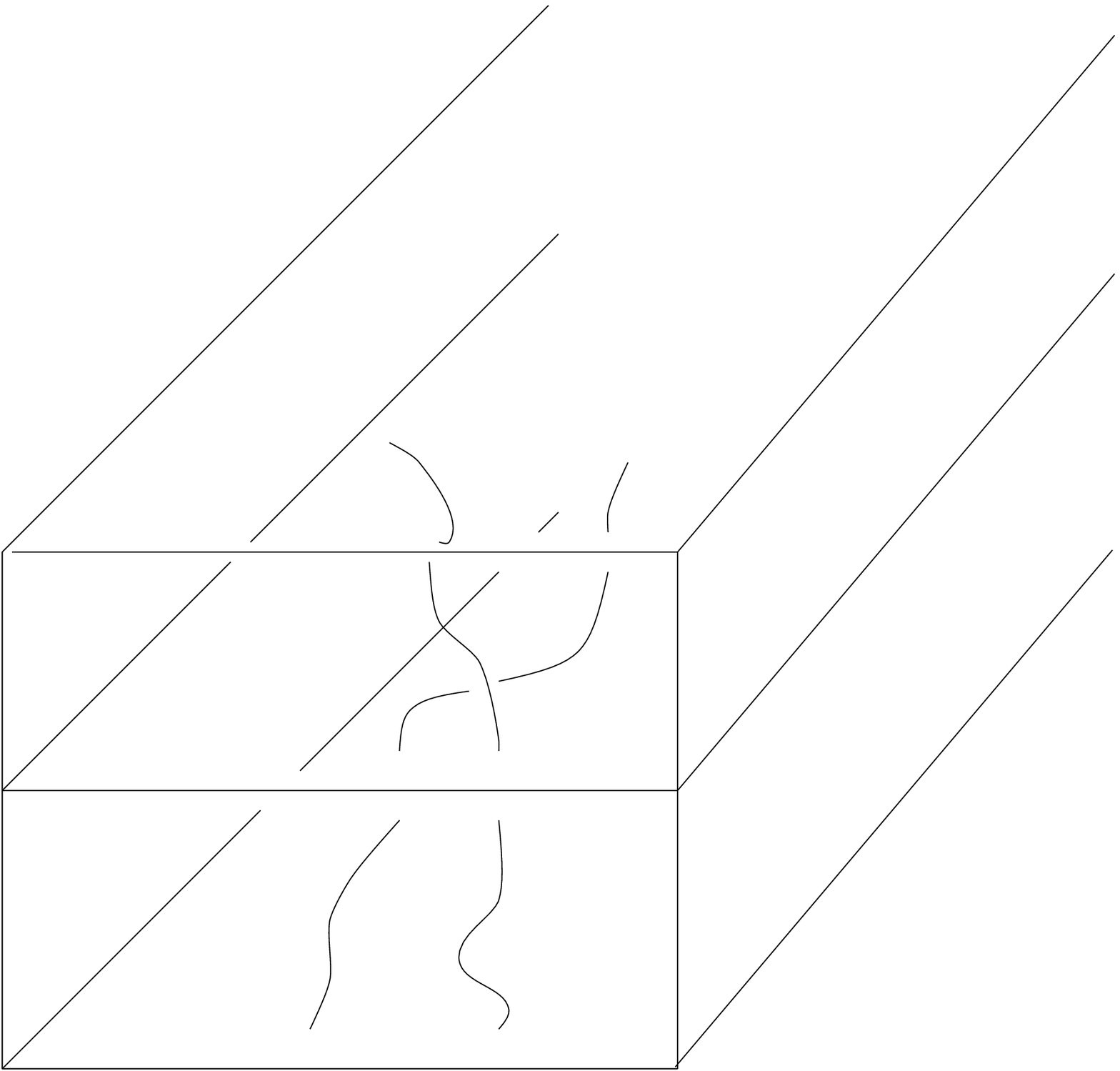}
\put(-260,0){$t=1$}
\put(-260,35){$t=0$}
\put(-260,70){$t=-1$}
\put(-150,-20){$\bdy$}
\caption{The situation for the proof of Lemma 2}
\end{figure}

Our goal is to show that we lose no generality in assuming that
$H$ is level-preserving in $\bdy \times [0,1]$.  This happens in
two steps.  First we make it preserve the $\bdy \times \{0\}$
level and then the $\bdy \times (0,1)$ levels.   We do this by
reparametrizing the arcs $L \cap (\bdy \times [0,1])$ one at a
time.  Denote by $t$ the elements of $[- \epsilon, \epsilon]$ and
by $s$ those of $[\delta, \eta]$

For the first step define a function $\tau : [\delta, \eta]
\rightarrow [- \epsilon, \epsilon]$ by \[\tau(s)=f_s^{-1}(f_s([-
\epsilon, \epsilon]) \cap (\bdy \times \{0\}))\] This is well-
defined since $h \circ f_s$ is increasing.  We want a function
$F:[- \epsilon, \epsilon] \times [\delta, \eta] \rightarrow [-
\epsilon, \epsilon]$ so that
\begin{enumerate}
 \item $F_s$ is a diffeomorphism for all $s$, and
 \item $F_s(\tau (s)) = 0$
\end{enumerate}  
Choose a small $\gamma >0$ so that $(\tau (s) - \gamma , \tau (s)
+ \gamma) \subset [-\epsilon, \epsilon]$ for each $s$. Let $\beta
_s : [- \epsilon, \epsilon] \rightarrow [0,1]$ be a smooth
function with $\beta
_s(\tau (s)) = 1$ and $\beta _s(t) = 0$ outside of $(\tau (s) -
\gamma , \tau (s) + \gamma)$.  Then define $F_s(t) = t - \tau (s)
\beta _s(t)$.  Clearly $F_s(\tau (s)) = 0$.  The proof that $F_s$
is a diffeomorphism is rather technical and is omitted.  The
desired parametrizations are given by $g_s = f_s \circ F_s^{-1}$.

Now we want to make $H$ preserve levels in $\bdy \times [0,1]$. 
To this end we reparametrize again so that each arc is
parametrized by height.  Fortunately this is much easier than the
first part. Notice that for $h \circ g_s : [- \epsilon , 0]
\rightarrow [0,1]$ we have
\begin{enumerate}
 \item $h \circ g_s(- \epsilon) = 0$
 \item $h \circ g_s(0) = 1$, and
 \item $(h \circ g_s)^{\prime}(t) \neq 0$ for all $t$
\end{enumerate}
Thus $h \circ g_s$ is a diffeomorphism for each $s$.  Let $\zeta
_s:[0,1] \rightarrow [- \epsilon ,0]$ be $(h \circ g_s)^{-1}$,
and consider $g_s \circ \zeta _s : [0,1] \rightarrow M$, a
parametrization of one of the arcs.  Since $h \circ g_s \circ
\zeta _s = 1_{[0,1]}$ the arc is parametrized by height for all
$s$.  Clearly we can do this for each arc in $\bdy \times [0,1]$,
so there is no harm in assuming that $H$ is level preserving in a
neighborhood of $\bdy$.

Next we define $G$. For ease of notation let's assume that
$\delta = 0$ and $\eta = 1$. Let $b:I \rightarrow I$ be a smooth 
function with $b(0)=1$ and $b(1)=0$.  Then make 
\[ G(x,t) = \left\{ \begin{array}{ll}

                     x & \mbox{$x \in M-(\partial M \times I),    
                          0 \leq t \leq \frac{1}{2}$} \\
                     H(x,2b_1(s)t) &\mbox{$x=(y,s) \in \partial M 
                                           \times I, 0 \leq t     
                                           \leq \frac{1}{2}$} \\
                     H(x, 2t-1) &  \mbox{$x \in M-(\partial M    
                                   \times I), \frac{1}{2}<t \leq  
                                   1$} \\
                     H(x, (2-2b_1(s))(t-1)+1) & \mbox{$x=(y,s)    
                                   \in \partial M \times I,       
                                       \frac{1}{2}<t \leq 1$}
                     \end{array}
            \right. \]
It is easy to check that $G(L,1)=H(L,1)$ and that the maps line
up properly at $t=\frac{1}{2}$ and at the boundary of $\partial M
\times I$ that is interior to $M$.  Notice that the only way that
$G$ differs from $H$ is in the speed at which the isotopy
happens.  This means that our assumption that $H$ is level
preserving assures us that $H(L,t)$ is actually an embedding for
all $t \in I$. Also it isn't hard to see that $G(L,\frac{1}{2}) =
\sigma \cdot L$ and $G$ fixes $\partial M$ for all $t \in
[\frac{1}{2}, 1]$. 
\end{proof}

Ying-Qing Wu suggested a much shorter proof of Lemma 2 but in the proof of the
main theorem we want to use this technique again, so we use the longer 
version.

\section{Proof of the Theorem}
\begin{proof}
The first step is to define an $R$-module
homomorphism from ${\cal L}_{n}(X_{1}) \times {\cal
L}_{n}(X_{2})$ to ${\cal L}(X)$.  If we just take the union the
framed arcs becomed framed circles.  Call this map $q_{n}$. It
is well-defined since the framed braids on either side are nicely
embedded. Extend this to a map \[q_{n}:R{\cal L}(X_{1}) \times
R{\cal L}(X_{2}) \rightarrow R{\cal L}(X)\] bilinearly.  This
induces a homomorphism \[\hat q_{n}:R{\cal L}(X_{1}) \otimes
R{\cal L}(X_{2}) \rightarrow R{\cal L}(X)\] Now let $S$ be the
submodule of $R{\cal L}(X_{1}) \otimes R{\cal L}(X_{2})
$ generated by all elements of the form $a^{\prime} \otimes b$ or
$a \otimes b^{\prime}$, where $a \in
R{\cal L}(X_{1})$, $b \in R{\cal L}(X_{2})$, $a^{\prime} \in
K^1_n$, and $b^{\prime} \in K^2_n$. (Recall that $K_n^i$ is the smallest
submodule of $R{\cal L}_n(X_i)$ containing all of the skein and
framing relations.)
Notice that $\hat q_{n}(S)
\subset K$.  Therefore there is a homomorphism \[ \tilde q_{n}:
[R{\cal L}(X_{1}) \otimes R{\cal L}(X_{2})]/S \rightarrow K(X)\]
But $(R{\cal L}(X_{1}) \otimes R{\cal L}(X_{2}))/S \cong K_n(X_1)
\otimes K_n(X_2)$.  So we get a map \[\tilde
q_n:K_n(X_1) \otimes K_n(X_2) \rightarrow K(X)\]
   Thus for each $n$ there are embeddings \[\varphi_n:[K_n(X_1)
\otimes K_n(X_2)]/ \mathrm{ker}\tilde q_n \rightarrow K(X)\]
This implies the existence of \[\Phi:\bigoplus_{n=0}^{\infty}
\left[ K_n(X_1) \otimes K_n(X_2) \right] / \mathrm{ker}\tilde
q_n \rightarrow K(X)\] 
By Lemma 1 $\Phi$ is onto  
\[ \Rightarrow \bigoplus_{n=0}^{\infty}
 \frac{[K_n(X_1) \otimes K_n(X_2)]}{\mathrm{ker} \tilde
q_n} / \mathrm{ker}\Phi \cong K(X)\]  

\[ \Rightarrow \frac{\bigoplus_{n=0}^{\infty}
[ K_n(X_1) \otimes K_n(X_2)]}{ \bigoplus_{n=0}^{\infty}
\mathrm{ker} \tilde q_n} / \mathrm{ker} \Phi \cong
K(X)\]

This is a start but the result will not be very useful until we
characterize $\bigoplus_{n=0}^{\infty} \mathrm{ker} \tilde
q_n$ and $\mathrm{ker} \Phi$ topologically.  The first will be
described with the help of the group action discussed earlier and
we will deal with the second via bigon moves.

We want to show that if $L \in {\cal L}(X)$ and $H:L \times I
\rightarrow X$ is an isotopy then $H_0(L)$ and $H_1(L)$ are
related by a sequence of braiding moves and bigon moves. There is no harm in 
assuming that $H_0(L) \cap F = \{a_1, \ldots ,a_{2n} \} = H_1(L) \cap F$.
Furthermore we reparameterize the arcs near their intersection with $F$
using the technique of the proof of Lemma 2.  This allows us to view $H$
as two separate isotopies; one on $X_1$ and one on $X_2$.
We also assume 
that the motion $H$ is generic.  This means that there exist
$s_1, \ldots ,s_k \in I$ for which $|H_{s_i}(L) \cap F|$
is odd, and further only one of these points is non-generic.

\begin{figure}[h]
\centering
\leavevmode
\epsfxsize = 3.5in
\epsfysize = 1in
\epsfbox{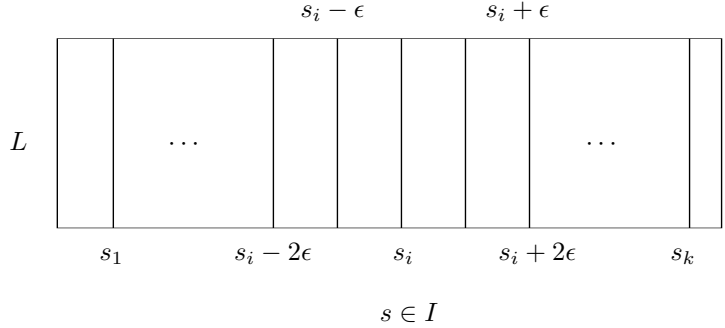}
\put(-130,-35){$s \in I$}
\put(-270,30){$L$}
\put(-125,-12){$s_i$}
\put(-160,80){$s_i-\epsilon$}
\put(-185,-12){$s_i-2\epsilon$}
\put(-236,-12){$s_1$}
\put(-20,-12){$s_k$}
\put(-85,-12){$s_i+2\epsilon$}
\put(-90,80){$s_i+\epsilon$}
\put(-210,30){$\cdots$}
\put(-52,30){$\cdots$}
\caption{The isotopy $H$}
\end{figure}

Choose $\epsilon$ so small that $H:L \times [s_i-\epsilon ,
s_i+\epsilon] \rightarrow X$ is a bigon move.  Of course this
isn't technically the case since the other components of $L$ move
some, but by choosing $\epsilon$ small we can ignore this.   The
strategy is to speed up $H$ so that what occurs from $s=0$ to
$s=s_0-\epsilon$ occurs from $s=0$ to $s=s_0-2\epsilon$.  Then we
"undo" $H$ near $F$ from $s=s_0-2\epsilon$ to $s=s_0-\epsilon$. 
At this point we apply Lemma 3.2, which gives a braiding move. 
Then from $s=s_0-\epsilon$ to $s=s_0+\epsilon$ we have a bigon
move.  Then between $s=s_0+\epsilon$ and $s=s_0+2\epsilon$ we put
the
link back to how it was at $s=s_0-\epsilon$.  Repeat this process
for all of the $s_i$ and we're done.

Let $u=\frac{s_0-\epsilon}{s_0-2\epsilon}$ and define $G:L \times
[0,s_0-2\epsilon] \rightarrow X$ by $F(x,s)=H(x,us)$.  Then
choose a small product neighborhood $F \times [-1,1]$ of $F$ in
$X$ with $F=F \times \{0\}$.  Let $g:F \times I \rightarrow F$ be
an isotopy of $F$ taking $H_{s_0-\epsilon}(L) \cap F$ to $H_0(L)
\cap F$. Next let $b:[-1,1] \rightarrow [0,1]$ be a smooth
function with $b(0)=1$ and $b(1)=b(-1)=0$.  Consider $G:(F \times
[-1,1]) \times [s_0-2\epsilon, s_0-\epsilon] \rightarrow X$ given
by \[G(x,t,s) = (g(x,b(t)s),t)\].  Then $G$
satisfies the hypotheses of Lemma 3.2 and so is a braiding move. 
Extend $G$ to $[0, s_0+\epsilon]$ by tacking on the appropriate
bigon move.  Further extend $G$ by taking \[G(x,t,s)=(g(x,1-
b(t)s),t)\] for $s \in [s_0+\epsilon, s_0+2\epsilon]$.  Finally
repeat
this process as many times as necessary to see that $H_1(L)$ as
$H_0(L)$ with a sequence of braiding moves and bigon moves
applied to it.
\end{proof}

\section{The Case When $F$ is Nonseparating}
We now assume tha $X$ and $F$ are as before, but $F$ is nonseparating.
In other words, $X=X_1/\psi$ where $\psi:F \rightarrow F$ is the
diffeomorphism and $X_1$ is $X$ cut along $F$. With some minor changes in 
the definitions the theorem carries over to this setting.  What's more,
the same proof works. After cutting $X$ along $F$, $X_1$ contains two copies
of $F$, say $F, F^{\prime} \subset \partial X_1$.  This gives us two sets of 
embedded arcs, $a_j \subset F$, and $a_j^{\prime} \subset F^{\prime}$.

\begin{definition}
A collection of framed arcs in $X_1$ is nicely embedded if the ends of
the framed arcs coincide with the $a_j$ and the $a_j^{\prime}$ in a one-
to-one manner.
\end{definition}

Now we discuss the necessary changes in the definition of braiding
moves.  Here we define both a left and a right action of $B_n(F)$ on
${\cal L}_n(X_1)$.  The left action is the same as in definition 4 and is
again written $\sigma \cdot L$.  The right action is exactly the same 
but it happens near $F^{\prime}$ and is written $L \cdot \sigma$. (Here
$L \in {\cal L}_n(X_1)$ and $\sigma \in B_n(F) \cong B_n(F^{\prime})$.)
For each $n$ let $U_n$ be the submodule of $K_n(X)$ generated by all the
elements of the form  $a-\sigma \cdot a \cdot {\sigma}^{-1}$, where
$a \in K_n(X_1)$ and $\sigma \in B_n(F)$.  Then let $U=\oplus_{n=0}^{\infty}
U_n$ and call this the submodule generated by all braiding moves.

The definition of a bigon move is the same as before.  Here we define $V$
to be the submodule of $K_n(X_1)$ generated by all elements of the form
$a-a^{\prime}$ where $a^{\prime}$ is the result of applying a bigon move to $a$.

\begin{theorem}
Let $X$ be an orientable 3-manifold, $F \subset X$ a closed, orientable
nonseparating surface, and $\psi:F \rightarrow F$ a diffeomorphism.  Let
$X_1$ be $X$ cut along $F$, so that $X=X_1/\psi$. Then
$$K(X) \cong \frac{\oplus_{n=0}^{\infty} K_n(X_1)}{U \vee V}$$
\end{theorem}

\begin{proof}
This proof differs from that of Theorem 1 only in the algebraic part.  Obviously
there is no tensor product here.  Nonetheless this part of the proof is the 
same in spirit.  In lieu of writing it out, we just say that the map $\Phi$ is 
defined from $\oplus_{n=0}^{\infty}K_n(X_1)$ instead of from $\oplus_{n=0}
^{\infty}  K_n(X_1) \otimes K_n(X_2)$. The rest of the proof is identical
since it all occurs near the surface $F$.
\end{proof}

\end{document}